\documentclass[10pt]{article}
\usepackage{cite,amsfonts,a4wide}
\usepackage{amsmath}
\usepackage{float}
\input amssym.def
\input amssym.tex
\usepackage{geometry}
\newcommand{\BB}[1]{\mbox{${\Bbb #1}$}}

\newcommand {\be}{\begin{equation}}
\newcommand {\e}{\end{equation}}
\newcommand {\bes}{\begin{displaymath}}
\newcommand {\es}{\end{displaymath}}
\newcommand {\nit}{\noindent}

\newcommand {\bit}{\bibitem}

\newcommand {\MF}{{\rm F}}

\newcommand {\MZC}{{\rm ZC}}
\newcommand {\MZR}{{\rm ZR}}
\newcommand {\MBR}{{\rm BR}}

\newcommand{\Mfib}{{\rm fib}}
\newcommand{\Modfib}{{\rm odfib}}
\newcommand{\Mevfib}{{\rm evfib}}
\newcommand{\MONA}{{\rm ONA}}
\newcommand {\MPP}{{\rm PP}}
\newcommand {\MOP}{{\rm OP}}

\begin{document}


\begin{center}
{\bf \Large On the structure of the complement $\overline{\Mfib}$
of the set $\Mfib$ \\
of fibbinary numbers in the set of positive natural numbers \\
\vspace{0.5cm}
A.J. Macfarlane$\;$}
\footnote{{\it e-mail}: 
a.j.macfarlane@damtp.cam.ac.uk}\\
\vskip .5cm
\begin{sl}
Centre for Mathematical Sciences, D.A.M.T.P.\\
Wilberforce Road, Cambridge CB3 0WA, UK
\vskip .5cm
\end{sl}
\end{center}

\begin{abstract}

The set $\Mfib$ of fibbinary numbers is defined via a bijection between the
set $\BB{N}$ of natural numbers and $\Mfib$. Since the elements of the set
$\Mfib$ do not exhaust $\BB{N}$, the structure of the complement
$\overline{\Mfib}$ of $\Mfib$ in $\BB{N}$, is of interest. An explicit
expression $\overline{\Mfib}=\bigcup_{k \geq 1} \Phi_k$ is obtained in
terms of certain well-defined sets $\Phi_k, \; k \geq 1$. The key to its
proof lies in first considering the odd numbers involved in this statement: a
general treatment, with full justification, of the binary representations of
the odd numbers is developed, and exploited  in showing the expression quoted
for $\overline{\Mfib}$ to be correct. The main results of the article can also
be viewed as providing partitions of the set of natural numbers, and also of
its subset of odd numbers, that follow from the introduction of
the set $\Mfib$, and of its subset of odd integers.

\end{abstract}
$\quad$

$\quad$

\nit {\bf Keywords}: 
Fibonacci number, Zeckendorf representation, fibbinary number, complement,
binary representation

$\quad$
 
\nit {\bf Mathematical Subject Classification}: 11A63, 11A67

\section{Introduction}

A recent paper \cite{FWB}, describes the bijection
$\mathcal{Z}$ between the
set $\BB{N}$ of positive natural numbers and the set $\Mfib$ of fibbinary
numbers, see \cite{FWB,lsw}, and $\underline{A003714}$ in \cite{oeis},
which underlies the definition
of the positive integers $\Mfib(n), \; n \geq 1$.
The bijection \cite[Eqs. (1),(29)]{FWB} in question
\begin{align} {\mathcal{Z}}: & \;\; \BB{N} \rightarrow \Mfib \nonumber \\
 & \;\; n  \mapsto \Mfib(n)
 \label{aa1}, \end{align} \nit
depends on the use of Zeckendorf's theorem \cite{zeck,lek,CS},
\cite[p. 125]{as} and the
Zeckendorf representation $z(n)$ of $n \in \BB{N}$. The paper 
\cite{FWB} contains more information and detail on these matters than section
two here provides.

The elements of the set $\Mfib$ do not exhaust $\BB{N}$, and the purpose of
this paper is to describe the complement $\overline{\Mfib}$
in $\BB{N}$ of $\Mfib$
\be\label{aa2} \overline{\Mfib}= \BB{N} - \Mfib. \e\nit
In particular a set theoretic expression for  $\overline{\Mfib}$  is offered
and proved to be correct. In terms of sets $\Phi_k,\; k \geq 1$, defined in
section three, this formula reads as
\be\label{aa3}  \overline{\Mfib} = \bigcup_{k=1}^\infty \Phi_k. \e\nit

To approach (\ref{aa3}), introduce the set $\Modfib$, $\underline{A022341}$ in
\cite{oeis}, \cite{lsw,FWB}, which is the subset of $\Mfib$ which contains
all of its odd elements, and the sets $\Psi_k,\; k \geq 1$, consisting
of the all of the odd elements of $\Phi_k$. The major part of the effort in
completing the proof of (\ref{aa3})
resides in proving corrrect the expression
\be\label{aa4}  \overline{\Modfib} = \bigcup_{k=1}^\infty \Psi_k.\e\nit
for  the complement $\overline{\Modfib}$
\be\label{aa5} \overline{\Modfib}= \BB{N}_{\rm odd} - \Modfib. \e\nit
of $\Mfib$ in the set of positive odd numbers $\BB{N}_{\rm odd}$,
$\underline{A005408}$ in \cite{oeis}. After this proof
is completed, that of (\ref{aa2}) follows easily.

The results (\ref{aa2}) and (\ref{aa3}) can also be recast in the form
\be\label{aa6} \BB{N}= \Mfib+\bigcup_{k=1}^\infty \Phi_k, \e\nit
that is, as a partition of the natural numbers.
Similarly (\ref{aa4}) and (\ref{aa5})
provide a partition of the odd numbers. 
 
Section two summarises well-known material about the sets $\Mfib$ and $\Modfib$
and the bijection $\mathcal{Z}$, as  needed for later sections.
Section three defines
the sets $\Phi_k, \; k \geq 0$, adopting the convenient notation $\Mfib=\Phi_0$,
and introduces tables that give valuable insight into the structure of these
sets. The elements of the sets $\Psi_k, \; k \geq 0,$ with $\Psi_0=\Modfib$,
can be read off the results quoted for the sets $\Phi_k$ or off the associated
tables. Tables 1 and 2  relate to $\Mfib$ and $\Phi_1$. The latter is typical
of the tables that can be drawn up for $\Phi_k$ for any $k \geq 1$. Explicit
expressions are given for $\Phi_k(n), \Psi_k(n), \;  k \geq 1$, in terms of
$\Phi_0(n)=\Mfib(n), \Psi_0(n)=\Modfib(n)$.

In section four, the first odd number array, $\MONA1$ is presented as Table 3.
Its first or $k=0$ column contains the integers of $\Modfib$. Its $k$-th column,
$k \geq 1$, contains the integers of $\Psi_k$. It will be seen that all odd
integers are present each one exactly once in the array $\MONA1$. To prove
this the vital step is to construct the second odd number array, $\MONA2$. 
This is obtained from the array $\MONA1$ by replacing each of its elements by
its unique binary representation, $\MBR$, written in an appropriate way.
Proof of (\ref{aa4}), given in detail, exploits the information so presented
explicitly in $\MONA2$.

The display in $\MONA2$ of the $\MBR$s of its odd number entries is of
significance in its own right. It can be seen to provide a $\MBR$ interpretation
of the formula (\ref{ff6}) for $\Psi_k(n)$, and hence also for the formula
(\ref{ff5}) for $\Phi_k(n)$. See section five.

\section{The Zeckendorf representation and the fibbinary numbers}
\subsection{The Zeckendorf representation}

Let $\MF$ denote the set of Fibonacci numbers
\begin{align} \MF & =\{ F_n, \quad n \geq 2 \} \nonumber \\
& = \{1,2,3,5,8,13,21,34,55, \dots \}. \label{bb1} \end{align} \nit

Zeckendorf's theorem \cite{zeck,lek,CS,as} states that every
positive integer can be expressed uniquely as a sum of
non-consecutive Fibonacci numbers $F_n \in \MF$.
This leads to the Zeckendorf representation -- notation $z(n)$ and
abbreviation ZR -- of the integers $n \in \BB{N}, \; n \geq 1$.
This  can be written uniquely in the form
\be \label{bb2} z(n)= \sum_{i=2}^r a_i F_i, \e\nit
where
\begin{align}  & a_i \in \{ 0,1 \},   \quad a_r=1, \nonumber \\
  & a_i\;a_{i+1}=0,  
\quad 2 \leq i <r, \label{bb3}  \end{align} \nit
and $r$ is the largest integer such that $F_r \leq n$. The entry $a_2$ is
referred to as the least significant bit, LSB, and (\ref{bb2}) as the
Zeckendorf condition, $\MZC$.

For $z(n)$, given  by (\ref{bb2}), write also
\be\label{bb4} z(n)=(a_r,a_{r-1}, \dots, a_2)_F = (n)_F, \e\nit
where, within $(n)_F$, $n$ denotes the
string
\be\label{bb5} a_r a_{r-1} \dots a_2, \e\nit
of ones and zeros obtained from (\ref{aa2}). Define also the length $l(n)$
of the $\MZR$ of $n$, that is, the total number of ones and zeros
in the string (\ref{bb5}).

If the integers $n \in \BB{N}$ are listed in ZR form, it is seen that they
separate into Fibonacci subsets ${\mathcal{F}}_k$ of cardinality $F_k$ such
that $l(n)=k$ for all $n \in {\mathcal{F}}_k$.

\subsection{The fibbinary numbers}

The definition of $z(n)$ of (\ref{bb2}) motivates
introduction of the set
\be\label{bb6} \Mfib = \{\Mfib(n), \quad n=1,2, \dots \}, \e\nit
$\underline{A003714}$ in \cite{oeis},\cite{lsw,FWB},
of fibbinary numbers, wherein the positive integer $\Mfib(n)$ is defined by
\be\label{bb7} \Mfib(n)=\sum_{i=2}^r a_i 2^{i-2}, \e\nit
and use is made of the same set of
coefficients $a_i, \quad 2 \leq i \leq r$,
as are used in (\ref{bb2}). Eq. (\ref{bb7}) is the binary representation,
$\MBR$, of $\Mfib(n)$. As for the integers $n \in \BB{N}$ in ZR form
(\ref{bb2}), the $\Mfib(n)$
separate into Fibonacci subsets $\mathcal{F}_k$ of cardinality $F_k$ with
$l(\Mfib(n))=l(n)$. Accordingly the listing of the elements of $\Mfib$
is presented as
\begin{align}
\Mfib  = & \{\Mfib(n), \quad n=1,2,\dots \} \nonumber \\
 =  & \{ 1;2;4,5;8,9,10;16,17,18,20,21=\Mfib(12);32,33,34,36,37,40,41,
 42=\Mfib(20); \nonumber \\
 & 64,65,66,68,69,72,73,74,80,81,82,84,85=\Mfib(33); \nonumber \\
 & 128, \dots , 170=\Mfib(54); 256,\dots, 341=\Mfib(88); \dots \}, \label{b11}
\end{align}\nit
using semi-colons to indicate separation of the
entries into Fibonacci subsets. The largest element of the
$\mathcal{F}_k, \;\; k=1,2,\dots $,
subset of $\Mfib$ is $\Mfib(F_{k+2}-1)$, as noted explicitly in (\ref{b11})
for $k=5, \dots, 9$.

 Writing (\ref{bb7}) as
\be\label{bb8} \Mfib(n)=(a_r,a_{r-1} \dots a_2)_2 =(n)_2, \e\nit
in analogy with (\ref{bb5}), see the same string (\ref{bb4}) used for $n$
in $(n)_F$ and $(n)_2$ but with distinct meanings. Respect of the Zeckendorf
condition is a strict requirement in both uses. 

It is clear that the map (\ref{aa1}) from  $\BB{N}$ to $\Mfib$
defines a bijection
which may be expressed now more explicitly as
\begin{align} {\mathcal{Z}}: & \;\; \BB{N} \rightarrow \Mfib \nonumber \\
 & \;\; n =(n)_F \mapsto \Mfib(n)=(n)_2 \label{bb12}. \end{align} \nit

\subsection{The set $\Mfib$ and its subset $\Modfib$}

Define next the important subset $\Modfib$ of $\Mfib$ of odd fibbinary numbers,
$\underline{A022341}$ in \cite{oeis}, \cite{lsw,FWB}. This subset contains
all the odd elements of $\Mfib$ each one exactly once:
\begin{align} \Modfib  = &  \{ \Modfib(n), \quad n=0,1,\dots \}
\nonumber \\
 = & \{1,5;9;17,21;33,37,41;65,69,73,81,85; \nonumber \\
 & 129,133,137,145,149,161,165,169; \nonumber \\
 & 257, \dots ,341; 513, \dots \}. \label{bb14}
\end{align}\nit
The even fibbinary numbers of $\Mevfib$, $\underline{A022342}$ in
\cite{oeis}, are not used here. The elements of the subsets (\ref{bb14}) are
given in terms of those of $\Mfib$ by
\begin{align} \Modfib(n) & = 4\Mfib(n)+1 =(n01)_2, \; n=0,1, \dots,
\label{bb15} \\
\Mevfib (n) & = 2\Mfib(n)=(n0)_2, \; n=1,2, \dots,\nonumber \end{align} \nit
where $\Mfib(0)=0$  is used so that (\ref{bb15}) reads correctly for $n=0$.
Note also that (\ref{bb15}) gives
\be\label{bb13} l(\Modfib(n))=l(\Mfib(n))+2=l(n)+2. \e\nit

\section{The sets $\Phi_k,\; k\geq 1$, and their subsets $\Psi_k, \; k \geq 1.$}

The definition of all sets $\Phi_k, \; k \geq 0$  of interest here depends upon
$\quad$

\nit $\qquad\qquad$ a) their initial element $\alpha_k$, and

\nit $\qquad\qquad$ b) the following key property:

\nit If $j \in \Phi_k, \; k \geq 0$, then so also does the even integer $2j$,
and the odd integer$4j+1$, but $4j-1 \not\in \Phi_k$.

Given $\alpha_k$, the elements of $\Phi_k$  can be obtained by applying the
key property. For the set $\Phi_0=\Mfib$, with $\Phi_0(1)=\alpha_0=1$,
(\ref{b11})
gives a listing of its elements.

To construct the sets $\Phi_k, \; k \geq 1$ define their initial elements
\be\label{bb16} \alpha_k=4k-1, \e\nit
odd integers not present in $\Mfib$, and apply property $b)$.
This gives rise to the following:
\begin{alignat}{2}
\Phi_k = & \{ \Phi_k(n), \; k \geq 1, \quad n \geq 1, \; \Phi_k(1)=\alpha_k. 
\}, & \nonumber \\
\Phi_1 = & \{3;6;12,13;24,25,26;48,49,50,52,53=\Phi_3(12); & \nonumber \\
 & 96,97,98,100,101,104,105,106=\Phi_3(20);  &  \nonumber \\
 & 192,193,194,196,197,200,201,202,208,209,210,212,213=\Phi_3(33);\dots \},
 & \label{tr1}\\
\Phi_2  =& \{7;14;28,29;56,57,58;112,113,114,116,117; & \nonumber \\
&  224,225,226,228,229,232,233,234;448,\dots, 469; \dots \dots \},
& \label{tr2}\\
\Phi_3 = & \{11;22;44,45;88,89,90;176,177,178,180,181; &\nonumber \\
 & 352,353,354,356,357,360,361,362; 706, \dots, 725; \dots \}, &\label{tr3}\\
\Phi_4 = & \{15;30;60,61;120,121,122;240,241,242,244,245; & \nonumber \\
 & 480,481,482,484,485,488,489,490;960,\dots,981;\dots \}, &
 \label{tr4}\\
\Phi_5 = & \{ 19;38;76,77;152,153,154;304,305,306,308,309;608,
\dots,618;1216, \dots,1237;\dots \},
 & \label{tr5}\\
\Phi_6  = & \{ 23;46,47;92,93;184,185,186;368,\dots,373;736,\dots,746;
1472,\dots,1493;
\dots \}. \ & \label{tr6} \end{alignat}

The sets $\Psi_k, \; k \geq 0$ are defined to be the sets of all odd elements
of $\Phi_k, \; k \geq 0$. For $k=0$, $\Psi_0=\Modfib$, with elements shown in
(\ref{bb14}). The elements of $\Psi_k, \; k \geq 1$, may be read
off (\ref{tr1} -- \ref{tr6}).
\begin{alignat}{2}
\rule[-1mm]{0mm}{5mm} \Psi_k  = & \{ \Psi_k(n), \; k \geq 1, \quad
n \geq 0, & \nonumber \\
\rule[-1mm]{0mm}{5mm}\Psi_1  = & \{ 3;13; 25; 49,53;97,101,105;
193,197,201,209,213;\dots \}, & \label{tr11}\\
\rule[-1mm]{0mm}{5mm} \Psi_2  = & \{ 7;29;57;113,117;225,229,233;
449,453,457,465,469; \dots \},& \label{tr12}\\
\rule[-1mm]{0mm}{5mm} \Psi_3  = & \{ 11;45;89;177,181;
353,357,361; 705,709,713,721,725;\dots \}, & \label{tr13}\\
\rule[-1mm]{0mm}{5mm} \Psi_4  = &  \{15;61;121;241,245;
481,485,489;961,965,969,979,981; \dots \}, & \label{tr14}\\
\rule[-1mm]{0mm}{5mm} \Psi_5 = & \{ 19;77;153;305,309;605,609,613;
1217,1221,1225,1233,1237; \dots \}, & \label{tr15}\\
\rule[-1mm]{0mm}{5mm} \Psi_6  = & \{ 23;93;185;369,373;
737,741,747; 1473,1477,1481,1489,1493; \dots \}. & \label{tr16} \end{alignat}

As in \cite{FWB}, the structures of sets such as the sets
$\Phi_k, \;k \geq 0$, are seen clearly by presenting their elements in table
form. The rules governing the construction of such tables are fully
specified in section four of \cite{FWB}. Table 1, copied from
\cite[Table 2]{FWB}, gives the fibbinary table for $\Mfib=\Phi_0$.

\vspace{0.5cm}
\begin{center}
\begin{tabular}{|r|r|r|r|r|r|r|r|r|r|r|r|r|}
\hline
1   & & & & &  & & & &    & & &  \\
2   & & & & &  & & & &    & & &    \\
4   & & & & &  & & & 5    & & & & \\
8   & & & & &  9 & & & 10 & & & &   \\
16  & & & 17 & & 18 & & & 20 & & & 21 &   \\
32  & & 33 & 34 & & 36 & & 37 & 40 & & 41 & 42 &  \\
\hline
\rule[-2mm]{0mm}{6mm} 0 & 1 & 2 & 3 & 4 & 5 & 6 & 7 &
8 & 9 & 10 & 11 & 12   \\
\hline
\end{tabular}
\end{center}
\centerline{Table 1: The fibbinary table. }
\vskip 0.5cm

Table 2 gives the corresponding table for $\Phi_1$. This is typical of the
tables for the $\Phi_k,\; k \geq 1$.

\vspace{0.5cm}
\begin{center}
\begin{tabular}{|r|r|r|r|r|r|r|r|r|r|r|r|r|}
\hline
3   & & & & &  & & & &    & & &  \\
6  & & & & &  & & & &    & & &    \\
12   & & & & &  & & & 13    & & & & \\
24   & & & & &  25 & & & 26 & & & &   \\
48  & & & 49 & & 50 & & & 52 & & & 53 &   \\
96  & & 97 & 98 & & 100 & & 101 & 104 & & 105 & 106 &  \\
\hline
\rule[-2mm]{0mm}{6mm} 0 & 1 & 2 & 3 & 4 & 5 & 6 & 7 &
8 & 9 & 10 & 11 & 12   \\
\hline
\end{tabular}
\end{center}
\centerline{Table 2: The table for the set $\Phi_1$.}
\vskip 0.5cm

Various features of such tables are of use in section four.
First, for all $k \geq 0$,
at the head of all columns of each table all the odd integers of
$\Psi_k$ are found. Second, below each such odd integer, call it $j$, in the
same column, its even multiples $2^r j, \; r \geq 1$, are found.
Thus, if the $\MBR$
of $j$ is written $(j)_2$, then these  multiples have $\MBR$ $(j0^r)_2$.
This is important: once mastery of the $\MBR$s of the odd numbers has been
achieved, the same for the even numbers, and hence all numbers, follows easily.
Third, the elements of the Fibonacci subsets $\mathcal{F}_i, \;\; i=1,
\dots6$, of the set displayed are aligned along the rows $i=1,\dots6$,
of the table.

Next, note that the structure of the set $\Phi_0=\Mfib$
can be expressed in the form
\begin{align} \Phi_0= & \sum_{s=0}^\infty \alpha_0(s), \; \alpha_0(1)=1,
\nonumber \\
\alpha_0(s)= & \Modfib(s) \{2^r, \, r=0,1,\dots, \} \label{ff1}.\end{align}\nit
Similarly, for $k \geq 1$,

\begin{align} \Phi_k= & \sum_{s=0}^\infty \alpha_k(s),\; \alpha_k(1)=4k-1,
\nonumber \\
\alpha_k(s)= & \Psi_k(s) \{2^r, \, r=0,1,\dots, \} \label{ff2}.
\end{align}\nit

Alongside the result (\ref{bb15})
\be\label{ff3} \Psi_0(n)=\Modfib(n)= 4\Mfib(n)+1= 4\Phi_0(n)+1 \e\nit
it can be checked that 
\be\label{ff4} \Psi_k(n)=4\Phi_k(n)+1,  \; k \geq 1, \; n \geq 1,\e\nit
is correct.
Although $\alpha_k=4k-1$ holds, by definition, (\ref{bb16}),
(\ref{ff4}) cannot be extended to the case $n=0$, since $\Phi_k(0), \; k
\geq 1$, is not defined. Setting $\Phi_k(0)=k -\frac{1}{2}$ makes no
sense, so that $\Phi_k(0), \; k \geq 1$, is 
left undefined. This minor awkwardness is caused by the use of
accepted (and natural) definitions of $\Mfib$ and $\Modfib$.

An explicit expression for $\Phi_k(n), \; k \geq 1$, is easy to derive. It reads
as
\be\label{ff5}  \Phi_k(n)=\Phi_0(n)+(2k-1) 2^{l(n)}, \e\nit
where $l(j)$ is the length of $j$, or the total number of ones and zeros
in the $\MZR$ of $j$.
The result $\Phi_0(n)=\Mfib(n)$ does not conform to
(\ref{ff5}). Also for $k \geq 1,\; n \geq 0,$
\be\label{ff6}  \Psi_k(n)=\Psi_0(n)+(2k-1)2^{l(\Psi_0(n))}, \e\nit
and, again, $\Psi_0(n)=\Modfib(n), \; n \geq 0$, does not
conform to (\ref{ff6}).
The result (\ref{ff6}) is very important, in full agreement with the
picture of the $\MBR$s of the odd numbers developed in section four, and
commented on in section five. To show that (\ref{ff6}) is consistent
with (\ref{ff5}), use (\ref{bb15}) and (\ref{bb13}). See also the examples
(\ref{ex4}) and (\ref{ex5}) in section five.

\section{The odd numbers arrays $\MONA1$ and $\MONA$}

Using the definitions, (\ref{tr11} -- \ref{tr16}), of the sets
$\Psi_k, k \geq 0$, define  the first odd number array $\MONA1$, Table 3.

\begin{center}
\begin{tabular}{|r||r|r|r|r|r|r|}
\hline
\rule[-2mm]{0mm}{6mm} $\Psi_0$ & $\Psi_1$ & $\Psi_2$ & $\Psi_3$ &
$\Psi_4$ & $\Psi_5$ & $\Psi_6$ \\
\hline
1 & 3 & 7 & 11 & 15 & 19 & 23    \\
5 & 13 & 29 & 45 & 61 & 77  & 93  \\
9 & 25 & 57 & 89 & 121 & 153 & 185 \\
17 & 49 & 113 & 177 & 241 & 305 & 369  \\
21 & 53 & 117 & 181 & 245 & 309 & 373  \\
33 & 97 & 225 & 353 & 481 & 609 & 737 \\
37 & 101 & 229 & 357 & 485 & 613 & 741 \\
41 & 105 & 233 & 361 & 489 & 617 & 745  \\
65 & 193 & 449 & 705 & 961 & 1217 & 1473 \\
69 & 197 & 453 & 709 & 965 & 1221 & 1477 \\
73 & 201 & 457 & 713 & 969 & 1225 & 1481 \\
81 & 209 & 465 & 721 & 977 & 1233 & 1489 \\
85 & 213 & 469 & 725 & 981 & 1237 & 1493 \\
\hline
\end{tabular}
\end{center}
\centerline{Table 3: The first odd number array, $\MONA1$.}
\vskip 0.5cm

The columns of $\MONA1$ contain the integers of the sets $\Psi_k,\; k \geq 0$.
The rows of $\MONA1$ are numbered $n=0,1, \dots$, with entries $\Psi_0(0)
=\Modfib(0)=1$ and $\Psi_k(0)=\alpha_k=4k-1, \; k \geq 1$, in the top row.

The second odd number array $\MONA2$ is obtained by replacing each element
of $\MONA1$ by its unique binary representation, $\MBR$, giving rise to
Table 4 for $\MONA2$.

\begin{center}
\begin{tabular}{|l||l|l|l|l|l|l|}
\hline
\rule[-2mm]{0mm}{6mm} $\Psi_0$ & $\Psi_1$ & $\Psi_2$ & $\Psi_3$ &
$\Psi_4$ & $\Psi_5$ & $\Psi_6$ \\
\hline
0:1 & 1:1 & 11:1 & 101:1 & 111:1 & 1001:1 & 1011:1 \\
0:101 & 1:101 & 11:101 & 101:101 & 111:101 & 1001:101 & 1011:101 \\
0:1001& 1:1001 & 11:1001 & 101:1001 & 111:1001 & 1001:1001 & 1011:1001 \\
0:10001 & 1:10001 & 11:10001 & 101:10001 & 111:10001 & 1001:10001 & 1011:10001\\
0:10101 & 1:10101 & 11:10101 & 101:10101 & 111:10101 & 1001:10101 & 1011:10101\\
0:1$0^41$ & 1:1$0^41$ & 11:1$0^41$ & 101:1$0^41$ & 111:1$0^41$ & 1001:1$0^41$
& 1011:1$0^41$ \\
\hline
\end{tabular}
\end{center}
\centerline{Table 4: The second odd number array, $\MONA2$},
\centerline{showing the binary representations of the elements of $\MONA1$.}
\vskip 0.5cm

It is sometimes convenient to refer to the portion of an array
which lies to the right of the
double lines of the arrays as an $\MONA$ array proper.

The elements of $\Psi_0=\Modfib$ in the first column of $\MONA2$ respect the
Zeckendorf condition, $\MZC$, of (\ref{bb3}). By the definitions of $\Mfib$ and
$\Modfib$, all the elements of $\Mfib$ which do so are therefore present in
$\Phi_0$.
The elements of the sets $\Psi_k,\; k \geq 1$, all violate the $\MZC$, and hence
belong to $\overline{\Modfib}$. The $\MONA$ arrays proper in 
Table 3 and 4 consist entirely of odd numbers which violate the $\MZC$.
It is to be
proved that all integers of $\BB{N}_3$  which do so are present.

Since all elements of $\MONA2$ proper violate the $\MZC$,
their $\MBR$s must contain at least one pair $11$ of consecutive
entries equal to one. In $\MONA2$, the rightmost $11$ pair of any $\MBR$
is written as $1:1$, with the colon separating that $\MBR$
into two parts. To the left
of the colon, find the principal part $\MPP$ of the $\MBR$, with its rightmost
digit equal to one, so that, viewed itself as the $\MBR$ of some number, the
$\MPP$ represents an odd number. To the right of the colon, find the
$\Modfib$ part $\MOP$. The placing of the colon
ensures that the
$\MOP$ part of any $\MBR$ respects the $\MZC$, and represents an element
of $\Modfib$.
Next, note two observations regarding $\MONA2$ proper. First, the $\MPP$ of the
$\Psi_k(n), \; n \geq 0$, are independent of $n$ for each $k \geq 0$, with
these $\MPP$ equal to $0$  for $\Psi_0(n)$, and equal to $(2k-1)$ for
$\Psi_k(n), \; k \geq 1$.
Second, the $\MOP$ of the $\Psi_k(n), \; k \geq 0$, are independent of $k$,
with $\MOP$ equal to $\Modfib(n)$.

A colon has been inserted to
the left of the $\Modfib$ entries of the $k=0$ columns. Their $\MBR$s
therefore have $\MPP=0$.

The stage is now set to approach the proof of (\ref{aa4}). All the elements
of Table 4 are odd, by construction of the sets defining its columns.
Since $\MBR$s are unique no  odd number appears twice in Table 4. Put
otherwise, with entries written in $(\MPP,\MOP)$ form, any two distinct entries
differ either in $\MPP$ or $\MOP$ parts or both. Also, by construction,
the $\MPP$ part  of any entry is either $0$ for the $\Modfib$
elements of the $k=0$ columns,
or else $(2k-1)$  in columns $k \geq 1$, so that all possible odd numbers
are found for the $\MPP$ parts of all the entries of columns $k \geq 1$.
Thus, with the elements of Table 4 written in $(\MPP,\MOP)$ form, each part
ranges over all the allowed possibilities. No other odd number can arise by
the method of construction followed here. It is therefore claimed that all odd
numbers are present in Table 4. Those of $\Modfib$ occur in the $k=0$ column,
the rest occupy the other columns, so that the truth of (\ref{aa4}) follows.
To complete the justification of this statement, suppose it is alleged that
some odd number $\nu$ is missing from Table 4.  Either $\nu$ respects the
$\MZC$, in which case $\nu \in \Mfib$ and is 
indeed present in the first column, or else $\nu \in \overline{\Modfib}$,
in which case
its unique $\MBR$ must contains at least  one $11$ pair.
Find the rightmost such pair,
and identify the $\MPP$ and $\MOP$ parts of $\nu$. But all possible $\MBR$s
of the form
$(\MPP,\MOP)$ have been identified. So $\nu$ must coincide with exactly
one of them, and the allegation is refuted.

It is now claimed that the truth of (\ref{aa3}) follows from that of
(\ref{aa4}). In the sets $\Phi_0=\Mfib$, (\ref{b11}),
and $\Phi_k,\; k \geq 1$, (\ref{tr11} -- \ref{tr16}), it is clear from Tables 
such as Tables 1 and 2 that all the even elements of such sets are
found by reading down the columns of the tables headed by the odd entries.
It follows that all the even multiples of all the odd numbers are
accounted for, that is, all the even numbers, and hence all $n \in \BB{N}$.
Further they, all the integers $n \in \BB{N}$,
separate into the
elements of $\Mfib$, which respect the $\MZC$, and the rest which do not,
and hence belong to $\overline{\Mfib}$. Thus
(\ref{aa3}), and hence (\ref{aa2}), follow.

\section{ A few comments}
In Table 4  it is shown that the odd numbers $\Psi_k(n),\; k \geq 1,\;
n \geq 0$, can be expressed uniquely in the from 
\begin{alignat}{2} (\MPP,\;\MOP) = & ((2k-1), \; \Modfib(n)), \quad k \geq 1,
 & \nonumber \\
 = & (0, \; \Modfib(n)), \quad k=0.  & \label{gg1} \end{alignat}
 \nit
This picture is in full agreement with (\ref{ff6}). The second term of
(\ref{ff6}), interpreted in binary, correponds to $(2k-1)$ in  binary followed
by a number of zeros equal to the length of $\Modfib(n)$. This means binary
addition to this of the first term $\Modfib(n)$ is trivial, producing results
of the form seen in the entries of Table 4.

\nit Some examples, first for $\Psi$ sets and (\ref{ff6}):
\begin{alignat}{2}
57= & \Psi_2(2)=(11:1001)_2 =(1001)_2+(11)_2 2^4 & \nonumber \\
 = & 9 + 3. 16 = \Modfib(2)+(2.2-1)2^{l(\Modfib(2))}. &
 \label{ex1} \\
481= & \Psi_4(5)=(111:100001)_2 =
(100001)_2+(111)_2 2^6  & \nonumber \\
 = &  33 + 7. 64 =\Modfib(5)+(2.4-1)2^{l(\Modfib(5))}.
 & \label{ex2} \end{alignat}

\nit Examples for elements of $\Phi$ sets and (\ref{ff5}):
\begin{alignat}{2} 114= & \Phi_2(10)=(11:10010)_2 =(10010)_2+(11)_2 2^5
& \nonumber \\
 = &  18 + 3. 32 = \Mfib(10)+(2.2-1)2^{l(\Mfib(10)}. & \label{ex3} \\
485= & \Phi_4(17)=(111:100101)_2 =
(100101)_2+(111)_2 2^6  & \nonumber \\
 = &  37 + 7. 64  =\Mfib(17)+(2.4-1)2^{l(17)} & \label{ex4} \end{alignat}\nit
Note that (\ref{ex4}) can also be viewed in the context of (\ref{ff6}):
\begin{alignat}{2}
485 = & \Psi_4(6) = 37+448 & \nonumber \\
= & \Modfib(6) +(2.4-1)2^{l(\Modfib(6)},  & \label{ex5} \end{alignat}
with $\Mfib(17)=37=\Modfib(6)$.

\end{document}